\documentclass[12pt]{article}

\usepackage{amsmath}                            
\usepackage{amssymb}                            
\usepackage{amsthm}   
\usepackage{adjustbox}   
\usepackage{xcolor}

\newcommand{\defn}[1]{\textbf{#1}}
\newcommand\sq{\mathbin{\text{\scalebox{.84}{$\square$}}}}

\newtheorem{theorem}{Theorem}
\newtheorem{conjecture}[theorem]{Conjecture}
\newtheorem{lemma}[theorem]{Lemma}
\newtheorem{claim}{Claim}

\title{All group-based latin squares possess near transversals}

\author{Luis Goddyn, Kevin Halasz\footnote{Department of Mathematics, Simon Fraser University, 8888 University Dr, Burnaby, BC V5A 1S6, Canada}}

\begin{document}

\maketitle

\begin{abstract}
In a latin square of order $n$, a near transversal is a collection of $n-1$ cells which intersects each row, column, and symbol class at most once. A longstanding conjecture of Brualdi, Ryser, and Stein asserts that every latin square possesses a near transversal. We show that this conjecture is true for every latin square that is main class equivalent to the Cayley table of a finite group.

\end{abstract}


\section{Introduction and main theorems}

A \defn{latin square} of order $n$ is an $n \times n$ array in which each row and column is a permutation of some set of $n$ symbols. We refer to the set of cells containing any fixed symbol as a \defn{symbol class}. Let $L = [L_{i,j}]$ be a latin square of order $n$. 
A \defn{partial transversal} of $L$ is a collection of cells which intersects each row, column, and symbol class at most once. A \defn{transversal} is a partial transversal of size $n$ and a \defn{near transversal} is a partial transversal of size $n-1$. Although it is straightforward to find latin squares possessing no transversals (see \cite[p.\! 405]{Wanless2011}), there is no known example of a latin square which does not possess a near transversal. 

\begin{conjecture}\label{Brualdi}
Every latin square possesses a near transversal.
\end{conjecture}

First discussed in the literature roughly 50 years ago, Conjecture \ref{Brualdi} has been variously attributed to Brualdi, Ryser, and Stein (see \cite[Section 5]{Wanless2011}). The strongest general lower bound to date is due to Hatami and Shor \cite{HatamiShor2008}, who showed that every latin square possesses a partial transversal of size $n-O\left(\log^2(n)\right)$. 
There have also been numerous attempts to establish this conjecture as a special case of some stronger statement, including work in terms of hypergraph matchings \cite{AharoniBerger2009}, covering radii of sets of permutations \cite{CameronWanless2005}, and colorings of strongly regular graphs \cite{Goddynetal2019}. The present paper approaches Conjecture \ref{Brualdi} from the opposite direction by proving its most widely discussed special case (see \cite[p.\! 335]{DenesKeedwellBook}).

 Let $G$ be a finite group. The \defn{Cayley table} of $G$, denoted $L(G)$, is the latin square with rows and columns indexed by the elements of $G$ where $L(G)_{g,h}=gh$. In general, a latin square is said to be \defn{group-based} if there is an interpretation of its symbols which allows it to be realized as the Cayley table of some finite group. 

\begin{theorem}\label{mainthm}
Every group-based latin square possesses a near transversal.
\end{theorem}

It is worth noting that Theorem \ref{mainthm} can be trivially extended to a slightly wider class of latin squares. The existence of a near transversal is not affected by relabelling the rows, columns, or symbols of $L$, nor is it affected by permuting the roles played by rows, columns, and symbols. Thus, every latin square which is \defn{main class equivalent} to a group-based latin square possesses a near transversal by Theorem \ref{mainthm}.

%
%
%
%
%
%
%
%
%
%
%

We prove Theorem \ref{mainthm} using a graph-theoretic technical lemma. Letting $\mathcal{I}$ denote the set indexing the rows and columns of the latin square $L$, the \defn{latin square graph} $\Gamma(L)$ is defined on the vertex set $\{(r,c) \,:\, r,c \in \mathcal{I}\}$ with $(r,c) \sim (s,d)$ if and only if one of $r=s$, $c=d$, or $L_{r,c} = L_{s,d}$ holds. Note that there is a natural tripartition of $\Gamma(L)$'s edges into, respectively, \defn{row edges}, \defn{column edges}, and \defn{symbol edges}.  
Moreover, there is a bijective correspondence between near transversals of $L$ and independent sets of size $n-1$ in $\Gamma(L)$.


Given graphs $\Gamma = (V,E)$ and $\Gamma^\prime=(V^\prime,E^\prime)$, the \defn{disjoint union} of $\Gamma$ and $\Gamma^\prime$ is 
$\Gamma+\Gamma^\prime:=(V \sqcup V^\prime, E \sqcup E^\prime).$
For a positive integer $k$, we write $k\Gamma$ for the disjoint union of $k$ copies of $\Gamma$. Given
a set $W \subseteq V$, the \defn{induced subgraph} of $\Gamma$ with respect to $W$ is 
$\Gamma[W] := (W,\{e \in E: e \subseteq W\}).$
The \defn{M\"{o}bius ladder} of order $2n$, denoted $M_n$, is the cubic graph formed from a cycle of length $2n$---referred to as the \emph{rim} of $M_n$---by adding $n$ edges, one joining each pair of vertices at distance $n$ in the initial cycle. 
The \defn{prism graph} of order $2n$, denoted $Y_n$, is the Cartesian product of a cycle of length $n$ and an edge; symbolically $Y_n := C_n \sq K_2$.

\begin{lemma}\label{heavylifting}
Let $L$ be a group-based latin square of even order $n$, let $k$ be the greatest power of 2 dividing $n$, and let $l:=n/k$. 
%
If $L$ does not possess a transversal, then there is a positive integer $m$ dividing $l$ such that $\Gamma(L)$ has an induced subgraph isomorphic to 
\[\Lambda_{n,m} := M_{km} + \left(\frac{l-m}{2}\right)Y_{2k} .\]
\end{lemma}

Section \ref{mostwork} will be devoted to proving Lemma \ref{heavylifting}, while Section \ref{conclusion} mentions two possible extensions of the present work. We conclude this section by proving Theorem \ref{mainthm} assuming Lemma \ref{heavylifting}.

\begin{proof}[Proof of Theorem \ref{mainthm}]
Let $L$ be a group-based latin square of order $n$. We may assume $L$ does not possess a transversal. As first shown in \cite{Bruck1944}, this implies $n$ is even. We may therefore apply Lemma \ref{heavylifting} to find an induced copy of $\Lambda_{n,m}$ in $\Gamma(L)$. 
Because the $\left(2k(l-m)\right)$-vertex graph $\left(\frac{l-m}{2}\right)Y_{2k}$ is bipartite, it contains an independent set of size at least $k(l-m)$. Moreover, one can find an independent set of size $km-1$ in $M_{km}$ by greedily selecting vertices in cyclic order around its rim. Thus $\Lambda_{n,m}$ contains an independent set of size $(l-m)k + km = n-1$ which corresponds to a near transversal of $L$.
\end{proof}

\section{Proof of Lemma \ref{heavylifting}}\label{mostwork}

Let $G$ be a group of order $n$ with identity element $e$ and let $Syl_2(G)$ denote the isomorphism class of $G$'s Sylow 2-subgroups. A bijection $\sigma: G \rightarrow G$ is called a \defn{complete mapping} of $G$ if the set $\{L(G)_{g,\sigma(g)} \,:\, g \in G\}$ is a transversal of $L(G)$. 
%
%
Using the classification of finite simple groups and partial results of Hall and Paige, Bray, Evans, and Wilcox recently characterized the groups possessing complete mappings. 
\begin{theorem}[\cite{Evans2009, HallPaige1955, Wilcox2009}]\label{EvansWilcox}
A finite group $G$ possesses a complete mapping if and only if $Syl_2(G)$ is either trivial or non-cyclic.
\end{theorem}

It is not hard to check that if $H$ is a group of odd order, then the identity map is a complete mapping. Another nicely structured complete mapping for groups of odd order was found by Beals, Gallian, Headley, and Jungreis.
\begin{lemma}[\cite{Bealsetal1991}]\label{harmonious}
For every group $H$ of odd order $m$, there exists an ordering $H = \{h_0,h_1,\ldots,h_{m-1}\}$ such that, taking indices modulo $m$, both $h_i \mapsto h_{i+1}$ and $h_i \mapsto h_i$ are complete mappings.
\end{lemma}

Given two subsets $X_1,X_2 \subseteq G$ the \defn{product set} of $X_1$ by $X_2$ is
\[X_1X_2 := \{x_1x_2\,:\, x_1 \in X_1, \, x_2 \in X_2\}.\]
%
We write $Xy$ for the product set $X\{y\}$. Let $K$ be a subgroup of $G$ and let $H$ be a normal subgroup of $G$. We say that $G$ is \defn{the semidirect product} of $K$ and $H$, written $G =  K \ltimes H$,  if $K \cap H = \{e\}$, $KH = G$, and $|G| = |K||H|$.
The following was noted in \cite{HallPaige1955} as following from a result of Burnside.
\begin{lemma}[Burnside, \cite{HallPaige1955}]\label{sdp}
Let $G$ be a finite group and let $K$ be a Sylow 2-subgroup of $G$. If $K$ is cyclic and nontrivial, then there is a normal subgroup of odd order $H \triangleleft G$ such that 
\[G=  K \ltimes H.\]
\end{lemma}
%
%


To simplify notation we set $[n] := \{0,1,\ldots,n-1\}$ for every positive integer $n$.

\begin{proof}[Proof of Lemma \ref{heavylifting}]
Let $L$ be a latin square based on a group $G$ of order $n=kl$, where $k\geq 2$ is a power of 2 and $l$ is odd. Moreover, suppose $L$ does not possess a transversal. Theorem \ref{EvansWilcox} tells us that $Syl_2(G) = \mathbb{Z}_k$. It then follows from Lemma \ref{sdp} that $G$ contains a normal subgroup $H$ of order $l$ and an element $b$ of order $k$ such that 
\[G = \langle b \rangle \ltimes H.\]

 Let $a:= b^{k/2}$. As $H \triangleleft G$ and $a$ has order 2, $H$ has an automorphism 
 \[\alpha: h \mapsto aha.\]
 Let
 \[H^\ast := \{ h \in H \,:\, \alpha(h) = h\}\]
  and observe that $H^\ast$ is a subgroup of $H$. Let $m:=|H^\ast|$. As $m$ divides $l$ and $l$ is odd, $m$ is odd. By Lemma \ref{harmonious}, there is an ordering $H^\ast=\{h_0,h_1,\ldots,h_{m-1}\}$ for which the map $h_i \mapsto h_{i+1}$  is a complete mapping. Here and throughout the rest of this proof, indices are taken modulo $m$.

Let $\Gamma:= \Gamma(L)$. Toward defining a set $W \subseteq V(\Gamma)$ which induces $\Lambda_{n,m}$, let
\begin{align*} 
T_1 &:= \left\{ \left(b^ih_i,h_ib^i\right) \,:\, i \in [km] \right\}, \\
 T_2 &:=  \left\{ \left(b^ih_i,h_{i+1}b^{i+1}\right) \,:\, i \in [km] \right\},  \text{ and } \\
 T \,&:=T_1 \cup T_2. 
\end{align*}
Furthermore, let $F := H \setminus H^\ast$ and let
\begin{align*}
U_1 &:= \left\{ \left(b^if,fb^i\right) \,:\, f \in F,\, i \in [k]\right\},\\
U_2 &:= \left\{ \left(b^if,fb^{i+1}\right) \,:\, f \in F,\, i \in [k]\right\}, \text{ and }  \\
U\,  &:= U_1 \cup U_2.
\end{align*}
Finally, let
\[ W:=T \cup U.\]
We show $\Gamma[W] \cong \Lambda_{n,m}$ via the following series of three claims. 

\begin{figure}[t]
\begin{center}
\newcommand{\overred}[1]{\textcolor{red}{\mathbf{\underline{#1}}}}
\newcommand{\underblue}[1]{\textcolor{blue}{{\mathbf{#1}}}}
\newcommand{\blank}[1]{\textcolor{white}{#1}}
\def\arraystretch{1.2}
\adjustbox{scale=0.68}{
\begin{tabular}
{|@{\hspace{2pt}}c@{\hspace{2pt}}|@{\hspace{2pt}}c@{\hspace{2pt}}|@{\hspace{2pt}}c@{\hspace{2pt}}|@{\hspace{2pt}}c@{\hspace{2pt}}|@{\hspace{2pt}}c@{\hspace{2pt}}|@{\hspace{2pt}}c@{\hspace{2pt}}||@{\hspace{2pt}}c@{\hspace{2pt}}|@{\hspace{2pt}}c@{\hspace{2pt}}|@{\hspace{2pt}}c@{\hspace{2pt}}|@{\hspace{2pt}}c@{\hspace{2pt}}|@{\hspace{2pt}}c@{\hspace{2pt}}|@{\hspace{2pt}}c@{\hspace{2pt}}||@{\hspace{2pt}}c@{\hspace{2pt}}|@{\hspace{2pt}}c@{\hspace{2pt}}|@{\hspace{2pt}}c@{\hspace{2pt}}|@{\hspace{2pt}}c@{\hspace{2pt}}|@{\hspace{2pt}}c@{\hspace{2pt}}|@{\hspace{2pt}}c@{\hspace{2pt}}|}
  \hline
$\overred{1}$&$\overred{bc}$&$c^2$&$b$&$c$&$bc^2$&$d$&$d^2$&$cd$&$cd^2$&$c^2d$&$c^2d^2$&$bd^2$&$bd$&$bcd^2$&$bcd$&$bc^2d^2$&$bc^2d$\\\hline
$bc$&$\overred{c^2}$&$\overred{b}$&$c$&$bc^2$&1&$bcd$&$bcd^2$&$bc^2d$&$bc^2d^2$&$bd$&$bd^2$&$cd^2$&$cd$&$c^2d^2$&$c^2d$&$d^2$&$d$\\\hline
$c^2$&$b$&$\overred{c}$&$\overred{bc^2}$&1&$bc$&$c^2d$&$c^2d^2$&$d$&$d^2$&$cd$&$cd^2$&$bc^2d^2$&$bc^2d$&$bd^2$&$bd$&$bcd^2$&$bcd$\\\hline
$b$&$c$&$bc^2$&$\overred{1}$&$\overred{bc}$&$c^2$&$bd$&$bd^2$&$bcd$&$bcd^2$&$bc^2d$&$bc^2d^2$&$d^2$&$d$&$cd^2$&$cd$&$c^2d^2$&$c^2d$\\\hline
$c$&$bc^2$&1&$bc$&$\overred{c^2}$&$\overred{b}$&$cd$&$cd^2$&$c^2d$&$c^2d^2$&$d$&$d^2$&$bcd^2$&$bcd$&$bc^2d^2$&$bc^2d$&$bd^2$&$bd$\\\hline
$\overred{bc^2}$&1&$bc$&$c^2$&$b$&$\overred{c}$&$bc^2d$&$bc^2d^2$&$bd$&$bd^2$&$bcd$&$bcd^2$&$c^2d^2$&$c^2d$&$d^2$&$d$&$cd^2$&$cd$\\\hline\hline
$d$&$bcd^2$&$c^2d$&$bd^2$&$cd$&$bc^2d^2$&$\underblue{d^2}$&1&$cd^2$&$c$&$c^2d^2$&$c^2$&$\underblue{bd}$&$b$&$bcd$&$bc$&$bc^2d$&$bc^2$\\\hline
$d^2$&$bcd$&$c^2d^2$&$bd$&$cd^2$&$bc^2d$&1&$\underblue{d}$&$c$&$cd$&$c^2$&$c^2d$&$b$&$\underblue{bd^2}$&$bc$&$bcd^2$&$bc^2$&$bc^2d^2$\\\hline
$cd$&$bc^2d^2$&$d$&$bcd^2$&$c^2d$&$bd^2$&$cd^2$&$c$&$\underblue{c^2d^2}$&$c^2$&$d^2$&1&$bcd$&$bc$&$\underblue{bc^2d}$&$bc^2$&$bd$&$b$\\\hline
$cd^2$&$bc^2d$&$d^2$&$bcd$&$c^2d^2$&$bd$&$c$&$cd$&$c^2$&$\underblue{c^2d}$&1&$d$&$bc$&$bcd^2$&$bc^2$&$\underblue{bc^2d^2}$&$b$&$bd^2$\\\hline
$c^2d$&$bd^2$&$cd$&$bc^2d^2$&$d$&$bcd^2$&$c^2d^2$&$c^2$&$d^2$&1&$\underblue{cd^2}$&$c$&$bc^2d$&$bc^2$&$bd$&$b$&$\underblue{bcd}$&$bc$\\\hline
$c^2d^2$&$bd$&$cd^2$&$bc^2d$&$d^2$&$bcd$&$c^2$&$c^2d$&1&$d$&$c$&$\underblue{cd}$&$bc^2$&$bc^2d^2$&$b$&$bd^2$&$bc$&$\underblue{bcd^2}$\\\hline\hline
$bd$&$cd^2$&$bc^2d$&$d^2$&$bcd$&$c^2d^2$&$\underblue{bd^2}$&$b$&$bcd^2$&$bc$&$bc^2d^2$&$bc^2$&$\underblue{d}$&1&$cd$&$c$&$c^2d$&$c^2$\\\hline
$bd^2$&$cd$&$bc^2d^2$&$d$&$bcd^2$&$c^2d$&$b$&$\underblue{bd}$&$bc$&$bcd$&$bc^2$&$bc^2d$&1&$\underblue{d^2}$&$c$&$cd^2$&$c^2$&$c^2d^2$\\\hline
$bcd$&$c^2d^2$&$bd$&$cd^2$&$bc^2d$&$d^2$&$bcd^2$&$bc$&$\underblue{bc^2d^2}$&$bc^2$&$bd^2$&$b$&$cd$&$c$&$\underblue{c^2d}$&$c^2$&$d$&1\\\hline
$bcd^2$&$c^2d$&$bd^2$&$cd$&$bc^2d^2$&$d$&$bc$&$bcd$&$bc^2$&$\underblue{bc^2d}$&$b$&$bd$&$c$&$cd^2$&$c^2$&$\underblue{c^2d^2}$&1&$d^2$\\\hline
$bc^2d$&$d^2$&$bcd$&$c^2d^2$&$bd$&$cd^2$&$bc^2d^2$&$bc^2$&$bd^2$&$b$&$\underblue{bcd^2}$&$bc$&$c^2d$&$c^2$&$d$&1&$\underblue{cd}$&$c$\\\hline
$bc^2d^2$&$d$&$bcd^2$&$c^2d$&$bd^2$&$cd$&$bc^2$&$bc^2d$&$b$&$bd$&$bc$&$\underblue{bcd}$&$c^2$&$c^2d^2$&1&$d^2$&$c$&$\underblue{cd^2}$\\\hline
\end{tabular}
}%
\caption{Cayley table of $S_3 \times \mathbb{Z}_3= \langle b,c,d \,|\, b^2=c^3=d^3 = 1, bc = cb, bd = d^2b\rangle$ with $\overred{T}$ and $\underblue{U}$ highlighted. Here $k=2$, $H=\langle c,d\rangle$, and $H^\ast = \langle c \rangle$, with $H^\ast$ ordered by $h_i = c^i$. The first six rows and columns are indexed by $\langle b \rangle H^\ast$ and the main diagonal is $T_1 \cup U_1$.}\label{figure}
\end{center}
\end{figure}

\begin{claim}\label{TandU}
$T \cap U =  \emptyset$ and there is no edge between $T$ and $U$.
\end{claim}

As $G = \langle b \rangle \ltimes H$, every element of $G$ has a unique representation of the form $g = b^ih$ for $i \in [k]$ and $h \in H$. Therefore, the definition of $F$ implies
\begin{equation} \label{Hpartition} \langle b \rangle H^\ast \cap \langle b \rangle F = \emptyset.\end{equation}
But for every $(t,s) \in T$ and every $(u,v) \in U$ we have $t \in \langle b \rangle H^\ast$ and $ u \in \langle b \rangle F$. Thus $T \cap U = \emptyset$ and there are no row edges between $T$ and $U$.

As $H \triangleleft G$ we have $b Hb^{-1} = H$. 
Moreover, as $a = b^{k/2}$, for every $h \in H^\ast$ we have $\alpha(bhb^{-1}) = bhb^{-1}$, so $bhb^{-1} \in H^\ast$. Thus
 \begin{equation}\label{normalizer}
b H^\ast b^{-1} = H^\ast.
 \end{equation}
It then follows from the definition of $F$ that
\begin{equation}\label{sgcommute}
bF b^{-1}= F
 \end{equation}
and, as the identity map is a complete mapping of both $H$ and $H^\ast$,
  \begin{equation}\label{sqcomp}
  f \mapsto f^2 \text{ is a permutation of  }F.
  \end{equation}
  Thus for every $(u,v) \in U$, both $v$ and $uv$ are in $\langle b \rangle F$. But \eqref{normalizer} tells us that for every $(t,s) \in T$, both $s$ and $ts$ are in $\langle b \rangle H^\ast$. It then follows from \eqref{Hpartition} that there are no column edges and no symbol edges between $T$ and $U$.

\begin{claim}\label{prisms}
$\Gamma[U]$ consists of $\frac{l-m}{2}$ disjoint copies of $Y_{2k}$. 
\end{claim}

Observe that, when enumerating the vertices in $U$, every element of $\langle b\rangle F$ occurs exactly twice as a first coordinate and, by \eqref{sgcommute}, exactly twice as a second coordinate. Thus, each vertex in $\Gamma[U]$ is incident to exactly one row edge and exactly one column edge, so that the row and column edges in $\Gamma[U]$ form a 2-factor (of $\Gamma[U]$). Specifically, they form $l-m$ disjoint $2k$-cycles $\{C_f\,:\, f \in F\}$, with each $C_f$ defined by the vertex-sequence
\begin{equation*}
(f,f), (f,fb), (bf,fb), \left(bf,fb^2\right), \ldots, \left(b^{k-1}f,fb^{k-1}\right), \left(b^{k-1}f,f\right).
\end{equation*}

It follows from the definitions of $H^\ast$ and $F$ that $\alpha|_F$ is a fixed-point free involution. Thus, to establish Claim \ref{prisms} it suffices to show that for every $i,j \in [k]$, every $f,h \in F$, and every $\epsilon,\delta \in \{0,1\}$, the vertices $(b^if,fb^{i+\epsilon})$ and $(b^jh,h b^{j+\delta})$ are joined by a symbol edge if and only if  $j \equiv i + k/2 \pmod{k}$, $h = \alpha(f)$, and $\epsilon = \delta$.

The ``if" direction of this equivalence follows directly from the definition of $\alpha$. For the converse direction we assume %
\[ b^if^2b^{i+\epsilon} = b^jh^2b^{j+\delta}\]
and, as latin square graphs are loopless, $(b^if,fb^{i+\epsilon}) \neq (b^jh,h b^{j+\delta})$.
It follows from \eqref{sgcommute} and \eqref{sqcomp} that
\[ b^if^2b^{i+\epsilon} \in b^{2i+\epsilon}F \text{ and } b^jh^2b^{j+\delta} \in b^{2j+\delta}F.\]
 Thus $\epsilon = \delta$ and $|i-j| \in \{0,k/2\}$. 
 
Now if $i=j$, then $b^if^2b^{i+\epsilon} = b^ih^2 b^{i+\epsilon}$ and \eqref{sqcomp} implies $f=h$, contradicting the fact that $(b^if,fb^{i+\epsilon})\neq(b^jh,h b^{j+\delta}) $. It follows that $j$ is the unique element of $[k]$ satisfying $j \equiv i+k/2 \pmod{k}$. Thus
\[b^ih^2b^{i+\epsilon} = b^{i+k/2} f^2 b^{i+\epsilon+k/2} = b^i \alpha(f^2) b^{i+\epsilon},\]
so $h^2 = \alpha(f^2) = (\alpha(f))^2$ and \eqref{sqcomp} implies $h = \alpha(f)$.

\begin{claim}\label{mobius}
$\Gamma[T]$ is isomorphic to $M_{km}$
\end{claim}

Observe that, when enumerating the vertices in $T$, every element of $\langle b\rangle H^\ast$ occurs exactly twice as a first coordinate and, by \eqref{normalizer}, exactly twice as a second coordinate. Thus, as is the case for $\Gamma[U]$, each vertex in $\Gamma[T]$ is incident to exactly one row edge and exactly one column edge. Unlike in $\Gamma[U]$, the row and column edges of $\Gamma[T]$ form a single cycle of length $2mk$. Indeed, as $m$ is odd and $|b|$ is a power of 2,
\[ (h_0,h_0), (h_0,h_1b),(bh_1,h_1b),\ldots, (b^{k-1}h_{m-1},b^{k-1}h_{m-1}), (b^{k-1}h_{m-1},h_0)  \]
is a Hamilton cycle in $\Gamma[T]$ which contains all of $\Gamma[T]$'s row and column edges. 

To establish Claim \ref{mobius} it suffices to show that for every $i,j \in [km]$ and every $\epsilon,\delta\in\{0,1\}$, the vertices $ (b^ih_i,h_{i+\epsilon}b^{i+\epsilon})$ and $(b^jh_j,h_{j+\delta}b^{j+\delta})$ are joined by a symbol edge  in $\Gamma[T]$  if and only if $i \equiv j + \frac{k}{2}m \pmod{km}$ and $\epsilon = \delta$. 

Indeed if $i \equiv  j + \frac{k}{2}m\pmod{km}$, then $i \equiv j \pmod{m}$ and, as $m$ is odd, $i \equiv j + k/2 \pmod{k}$. Together with $\epsilon  = \delta$,  as well as Lemma \ref{harmonious} and the definition of $H^\ast$, this implies
\[ b^ih_ih_{i+\epsilon}b^{i+\epsilon} = b^jah_jh_{j+\delta}ab^{j+\delta} = b^jh_jh_{j+\delta}b^{j+\delta},\]
which establishes the ``if" direction of the desired equivalence. 

For the converse direction consider $ (b^ih_i,h_{i+\epsilon}b^{i+\epsilon}), (b^jh_j,h_{j+\delta}b^{j+\delta}) \in T$  and assume that the group elements defining this pair of distinct vertices satisfy
\[ b^ih_ih_{i+\epsilon}b^{i+\epsilon} = b^jh_jh_{j+\delta}b^{j+\delta}.\]
From \eqref{normalizer} we see that
\[b^ih_ih_{i+\epsilon}b^{i+\epsilon} \in b^{2i+\epsilon}H^\ast \text{ and } b^jh_jh_{j+\delta}b^{j+\delta} \in b^{2j+\delta}H^\ast.\]
Thus $\epsilon = \delta$ and $i \equiv j \pmod{k/2}$. Now $b^j \in \{b^i,b^{i+k/2}\}$ and as $H^\ast$ is pointwise fixed by the automorphism $\alpha: h \mapsto b^{k/2} h b^{k/2}$,  both possible values of $b^j$ yield $h_ih_{i+\epsilon} = h_jh_{j+\epsilon}$. Lemma \ref{harmonious} then implies $h_i = h_j$, so $i \equiv j \pmod{m}$. 

Suppose $b^j = b^i$, which is equivalent to $i \equiv j \pmod{k}$. As $\gcd(k,m) = 1$ and $i,j \in [km]$, this implies $i=j$, contradicting the fact that $(b^ih_i,h_{i+\epsilon}b^{i+\epsilon})$ and $(b^jh_j,h_{j+\delta}b^{j+\delta})$ are distinct vertices. Therefore $j \equiv i+k/2 \pmod{k}$ and, as $k/2$ and $m$ are coprime, we conclude that $i \equiv j + \frac{k}{2}m \pmod{km}$. 
\end{proof}

\section{Concluding remarks}\label{conclusion}

The most obvious extension of the present paper is the general case of Conjecture \ref{Brualdi}. However, a proof of this conjecture would likely differ substantially from the argument presented above. Indeed, latin square graphs are in general not vertex-transitive, calling into question whether general latin square graphs can be shown to possess the sort of ``nice'' induced subgraphs found in Lemma \ref{heavylifting}.

There is, however, an extension of Theorem \ref{mainthm} to which the above techniques may be applicable. A partial transversal is \defn{non-extendable} if it is not contained in any larger partial transversal. The following conjecture was noted by Evans \cite[p.\! 470]{EvansBook} as a special case of a conjecture of Keedwell concerning sequenceable groups.
\begin{conjecture}[Keedwell]\label{keedwell}
For every finite non-Abelian group $G$, the latin square $L(G)$ possesses a non-extendable near transversal.
\end{conjecture}
We have established Conjecture \ref{keedwell} for those non-Abelian groups whose Sylow 2-subgroups are nontrivial and cyclic (these groups do not possess transversals, so a near transversal must be non-extendable). Perhaps our techniques can be used to find maximal independent sets of size $n-1$ in latin square graphs based upon non-Abelian groups with non-cyclic or trivial Sylow 2-subgroups. As far as we know, Conjecture \ref{keedwell} has not been attacked directly. However, many partial results are known due to its connection to sequenceable groups (see e.g. \cite{OllisSurvey}).





\bibliographystyle{plain} 
\bibliography{NTFJCTA}

 \vspace*{5mm}
 
\noindent
{\tt khalasz@sfu.ca, goddyn@sfu.ca}

\end{document}